# RELAXATION TIME OF *L*-REVERSAL CHAINS AND OTHER CHROMOSOME SHUFFLES

By N. Cancrini, P. Caputo and F. Martinelli

*Università di L'Aquila and INFM Roma, Università di Roma Tre and Università di Roma Tre*

We prove tight bounds on the relaxation time of the so-called *L*-reversal chain, which was introduced by R. Durrett as a stochastic model for the evolution of chromosome chains. The process is described as follows. We have $n$ distinct letters on the vertices of the $n$-cycle ($\mathbb{Z}$ mod $n$); at each step, a connected subset of the graph is chosen uniformly at random among all those of length at most $L$, and the current permutation is shuffled by reversing the order of the letters over that subset. We show that the relaxation time $\tau(n, L)$, defined as the inverse of the spectral gap of the associated Markov generator, satisfies $\tau(n, L) = O(n \vee \frac{n^3}{L^3})$. Our results can be interpreted as strong evidence for a conjecture of R. Durrett predicting a similar behavior for the mixing time of the chain.

**1. Introduction, models and results.** In a series of recent papers, R. Durrett has proposed stochastic models based on shuffling rules for the analysis of the evolution of chromosomes [6, 7, 8]. One of the main issues is that of determining the time needed for such processes to reach stationarity. Besides proving several results in this direction and discussing possible applications, he pointed out a number of interesting conjectures. One of them refers to the so-called *L*-reversal model, which is described as follows:

We consider a chain (the chromosome) of $n$ distinct genes, so that the configuration of the system at any given time is a permutation $\eta$ of $n$ letters over $n$ vertices. The vertices are lying on a circle and we take the graph structure of the $n$-cycle (where each vertex is connected to exactly two vertices). Given an integer $L \leq n$, one step of the *L*-reversal process is described as follows. Uniformly at random we pick a vertex $x$ and, independently, choose a number $1 \leq \ell \leq L$. Given $x$ and $\ell$, we perform the transition $\eta \to \eta^{x,\ell}$, where $\eta^{x,\ell}$ is obtained from $\eta$ by reversing the order of the letters in the segment $\{x, x+1, \ldots, x+\ell\}$, where the sums are taken modulo $n$. Note that the probability of the transition $\eta \to \eta^{x,\ell}$ coincides with that of the transition $\eta^{x,\ell} \to \eta$, and, therefore, we have reversibility with respect to the uniform distribution $\nu$ over all $n!$ configurations. As in [6], we consider the continuous-time version of the $L$-reversal dynamics, in which the steps described above are performed at the arrival times of an independent rate 1 Poisson process. This defines a continuous-time Markov chain converging to the uniform distribution $\nu$.

Let $T(n, L)$ denote the *mixing time* of the chain, that is,

$$T(n,L) = \inf\{t > 0 : \|p_t - \nu\|_{\mathrm{TV}} \leq 1/e\}, \tag{1}$$

where $\|\cdot\|_{\mathrm{TV}}$ stands for the total variation distance, and $p_t$ denotes the distribution of the configuration at time $t$ for the process started from some arbitrary (i.e., completely ordered) configuration $\eta$ (note that $\|p_t - \nu\|_{\mathrm{TV}}$ does not depend on the chosen $\eta$, because $\nu$ is the uniform measure).

The striking conjecture in [6] (see also Problem 4.1 in [8]) is that the mixing time $T(n, L)$ of the $L$-reversal process should satisfy

$$\frac{1}{C}\left(n \vee \frac{n^3}{L^3}\right) \log n \leq T(n, L) \leq C\left(n \vee \frac{n^3}{L^3}\right) \log n, \tag{2}$$

for some universal constant $C < \infty$. Here, for numbers $a, b$, we use $a \vee b$ to denote the maximum $\max\{a, b\}$. Note that, when $L = 1$, we have the usual transposition dynamics on the $n$-cycle, which has been studied in great detail; see, for example, [4, 5, 9, 10, 11]. In particular, the conjecture is known to hold in this case. Of most interest are the cases $L = n^\alpha$, $\alpha \in (0, 1)$.

To support the conjectured bounds in (2), Durrett uses an adaptation of Wilson's method to estimate mixing times [11]. In this way, by a careful choice of the slow modes of relaxation, he proves the lower estimates in (2). As for the upper estimates, by comparison with random transposition, he proves that

$$T(n, L) \leq C \frac{n^3}{L^2} \log n. \tag{3}$$

When $L = n^\alpha$, $\alpha \in (0, 1)$, the estimate (3) is off by a factor of order $n^\alpha$ for $\alpha \leq 2/3$ and $n^{2(1-\alpha)}$ for $\alpha \geq 2/3$. This shows that the comparison with transpositions is efficient only in the two extreme regimes $\alpha \in \{0, 1\}$. We are not able to prove (or disprove) the desired upper bound in (2), but we have considered the problem of establishing similar estimates for the *relaxation time* $\tau(n, L)$, that is, the inverse of the spectral gap of the chain. What we find confirms the predicted behavior, at least at the level of the spectral gap.

To describe our main result, we introduce some notation. For any function $f$ of the configurations, we write $\nu[f] = \frac{1}{n!}\sum_\eta f(\eta)$ for its expectation and $\mathrm{Var}(f) = \nu[f^2] - \nu[f]^2$ for its variance w.r.t. the uniform measure $\nu$. Also, we write $f^{x,\ell}$ for the function $\eta \to f(\eta^{x,\ell})$ and $R_{x,\ell}f = f^{x,\ell} - f$ for the associated gradient. The Dirichlet form of the $L$-reversal dynamics introduced above is

$$\mathcal{E}(f,f) = \frac{1}{2}\frac{1}{nL}\sum_{x=1}^{n}\sum_{\ell=1}^{L}\nu[(R_{x,\ell}f)^2]. \tag{4}$$

The relaxation time is then given by

$$\tau(n,L) = \sup_{f} \frac{\text{Var}(f)}{\mathcal{E}(f,f)}, \tag{5}$$

with the supremum taken over all nonconstant functions $f$.

Our main result reads as follows:

THEOREM 1.1. *There exists $C < \infty$ such that, for any $n \geq 2$ and any $1 \leq L \leq n$,*

$$\frac{1}{C}\left(n \vee \frac{n^3}{L^3}\right) \leq \tau(n,L) \leq C\left(n \vee \frac{n^3}{L^3}\right). \tag{6}$$

As usual, the lower estimate in (6) can be obtained by a suitable choice of test functions; see Section 3 for details. To obtain a sharp upper bound on $\tau(n,L)$, the argument is more delicate. Here we combine new comparison arguments with an improved iteration scheme inspired by recent work of Carlen, Carvalho and Loss [3]. We refer to the beginning of next section for a brief description of the main idea of the proof. Our method allows the derivation of sharp estimates for more general models in the class of so-called $p$-reversal processes, where each inversion $\eta \to \eta^{x,\ell}$ is performed with probability $p(\ell)$, $p$ being a probability vector on $\{1,\ldots,n\}$. The $L$-reversal process (4) is then obtained for $p(\ell) = \frac{1}{L}$, $\ell = 1,\ldots,L$. The so-called $\theta$-reversal process is defined by $p(\ell) \propto \theta^\ell$, $\theta \in (0,1)$, with Dirichlet form

$$\mathcal{E}_\theta(f,f) = \frac{1}{2}\frac{(1-\theta)}{n}\sum_{x=1}^{n}\sum_{\ell=1}^{n}\theta^{\ell-1}\nu[(R_{x,\ell}f)^2]. \tag{7}$$

In this case, calling $\tau(n,\theta)$ the relaxation time associated to (7) we obtain

THEOREM 1.2. *There exists $C < \infty$ such that, for any $n \geq 2$ and any $\theta \in (0,1)$,*

$$\frac{1}{C}(n \vee [n(1-\theta)]^3) \leq \tau(n,\theta) \leq C(n \vee [n(1-\theta)]^3). \tag{8}$$

We conclude this introduction with a final remark on Durrett's conjecture (2). It is well known that $T(n,L) \geq \tau(n,L)/C$. On the other hand, an application of standard bounds (see, e.g., Theorem 2.1.7 in [10]) only gives

$$T(n,L) \leq C\tau(n,L)\log(n!).$$

Better estimates can be obtained if one considers the decay of relative entropy functionals rather than $L^2$-norms. Namely, let $\beta(n,L)$ denote the *entropy dissipation constant* given by

$$\beta(n,L) = \sup_{f>0} \frac{\operatorname{Ent}(f)}{\mathcal{E}(f,\log f)}, \tag{9}$$

with the supremum taken over all positive nonconstant functions. Here we use the standard notation $\operatorname{Ent}(f) = \nu[f\log f] - \nu[f]\log\nu[f]$ for the entropy of $f$. It is not hard (see, e.g., Corollary 2.8 in [1]) to obtain the bound

$$T(n,L) \leq C\beta(n,L)\log\log(n!).$$

Therefore the conjecture (2) would follow if we could prove the bounds (6) with $\tau$ replaced by $\beta$. Despite recent progress in the understanding of entropy dissipation in various reversible Markov chains, the extension of our techniques to this case remains a challenging problem.

It is perhaps interesting to observe that the predicted mixing time cannot be derived from logarithmic Sobolev inequalities here. Indeed, we will show that if $s(n,L)$ denotes the log-Sobolev constant [obtained as in (9) with $\mathcal{E}(\sqrt{f},\sqrt{f})$ replacing $\mathcal{E}(f,\log f)$], then

$$s(n,L) \geq \frac{n^2}{CL}. \tag{10}$$

Comparing this with the bounds in Theorem 1.1, for $L = n^\alpha$, the bound (10) shows that $s(n,L) \geq \tau(n,L)n^{2\alpha-1}$ when $1/2 \leq \alpha \leq 2/3$ and $s(n,L) \geq \tau(n,L)n^{1-\alpha}$ when $\alpha \geq 2/3$. While comparison with random transposition can always be used to prove an upper bound of the form $s(n,L) \leq C(n^3/L^2)\log n$ (see, e.g., Theorem 6 in [6]), it remains an open question to understand the true asymptotic behavior in $L,n$ of the log-Sobolev constant.

The rest of this paper consists of two sections. In the first, we briefly outline the proof of the upper bound on the relaxation time, develop all the technical parts of our approach and finally prove the corresponding estimates of Theorem 1.1 and Theorem 1.2. In the final section, we prove the matching lower bounds and derive the estimate (10) on the log-Sobolev constant.

**2. Upper bounds.** The main idea needed for the proof of the upper bounds in Theorems 1.1 and 1.2 is to compare the $L$-reversal chain with (an average over $\ell \leq L$ of) a suitable "block dynamics" with blocks of size $\ell$ that we call the *block-average dynamics*. Such an auxiliary process seems to be at the heart of several shuffling chains, and it can be described informally as follows (see Section 2.3 below for the formal definition):

Assume for simplicity that $n = N\ell$ for some integers $N,\ell$ and partition the vertices $\{1,2,\ldots,n\}$ into $N$ nonoverlapping blocks $\{I_i\}$ of length $\ell$. Then, with rate $1/N$, a pair of blocks is chosen, say $I_i$ and $I_j$, and the letters inside $I_i \cup I_j$ are uniformly reshuffled over the vertices of $I_i \cup I_j$. We refer to this move as an "average" of the two blocks $I_i$ and $I_j$. A first key step then consists in proving that the relaxation time of the above block dynamics is $O(1)$ uniformly in $N, \ell$. We are actually able to compute exactly the spectral gap of this Markov chain (see Proposition 2.3 below). This computation is carried out via an adaptation and extension of the technique introduced by Carlen, Carvalho and Loss [3] for so-called Kac systems.

The comparison of the above-described block-average dynamics to the $L$-reversal chain is then accomplished in two further steps:

In the first step, we compare the block-average dynamics to an intermediate block dynamics in which only adjacent blocks are averaged but we also allow an "exchange" between two adjacent blocks. We say that blocks $I_i, I_{i+1}$ are *exchanged* if the configuration of letters inside $I_i$ is interchanged with that inside $I_{i+1}$. We shall call this intermediate dynamics the *local average-exchange block dynamics*. This step of the comparison uses rather standard path techniques, but the crucial point is that the usual diffusive scaling factor $N^2$ appears only in front of the exchange operations and does not multiply the local average moves (see Lemma 2.5).

In the second step, we compare the local average-exchange block dynamics to the $L$-reversal chain. The key to this step is the observation that moves of the $L$-reversal process can efficiently simulate both block exchanges and block averages. In fact, the local exchange move is easily expressed in terms of (three) suitable $\ell$-reversal moves $R_{x,\ell}$ (see Lemma 2.6). Also, the average of a pair of adjacent blocks $I_i$ and $I_{i+1}$ can be compared to inversions, using the Poincaré inequality for random transpositions on the complete graph generated by vertices in $I_i \cup I_{i+1}$. Finally, we use a straightforward comparison between the latter and a $2\ell$-reversal chain on $I_i \cup I_{i+1}$ (see Lemma 2.7). This does not spoil the estimate because, as discussed in the Introduction, the random transpositions chain and the $L$-reversal chain can be efficiently compared in the extreme case $L = n$ (i.e., $\alpha = 1$).

2.1. *Setting and notation.* The setting and frequently-used notation are described as follows. $V_n$ is the set of $n$ ordered vertices $\{1, 2, \ldots, n\}$. Unless otherwise stated, the sum $x + y$ is assumed to be taken modulo $n$, for any $x, y \in V_n$. A permutation (often called a *configuration*) of $n$ letters over $V_n$ is denoted by $\eta$, with $\eta_x = j$ meaning that we have letter $j$ at vertex $x \in V_n$. We write $\Omega$ for the space of all $n!$ configurations. $\nu$ is the uniform distribution on $\Omega$. For any integer $m \leq n$, $I \subset V_n$ is called an $m$-block if $I$ is of the form $\{x+1, x+2, \ldots, x+m\}$ for some $x \in V_n$. Given $k \in \mathbb{N}$, we write $\eta_{k,m}$ for the configuration over the $m$-block $\{(k-1)m+1, (k-1)m+2, \ldots, km\}$, that is, $\eta_{1,m}$ gives the letters over the first $m$-block, $\eta_{2,m}$ those over the second $m$-block and so on. Also, let $\Omega_m$ denote the set of possible configurations over the the first $m$-block (i.e., all possible realizations of $\eta_{1,m}$).

2.2. *Preliminary spectral estimates.* We start with a spectral computation that will be used in the arguments below. For any $m \le n$, let $\mathcal{K}$ denote the symmetric stochastic matrix given by the conditional probabilities

(11) $$\mathcal{K}(\alpha, \beta) := \nu[\eta_{1,m} = \beta | \eta_{2,m} = \alpha], \qquad \alpha, \beta \in \Omega_m.$$

In words, $\mathcal{K}(\alpha, \beta)$ stands for the probability of the configuration $\beta$ occuring in the first $m$-block, given that the configuration over the second $m$-block equals $\alpha$.

LEMMA 2.1. *For any $m \le n/2$, the spectrum of $\mathcal{K}$ is given by*

(12) $$\lambda_k := (-1)^k \binom{n-m-k}{m-k} \binom{n-m}{m}^{-1}, \qquad k = 0, \ldots, m.$$

PROOF. Let $\mathcal{H}_m$ denote the space of all functions $\varphi : \Omega_m \to \mathbb{R}$ and write

(13) $$\langle \varphi, \psi \rangle = \frac{1}{n!} \sum_{\eta \in \Omega} \varphi(\eta_{1,m}) \psi(\eta_{1,m}) = \frac{(n-m)!}{n!} \sum_{\alpha \in \Omega_m} \varphi(\alpha) \psi(\alpha)$$

for the scalar product of $\varphi, \psi \in \mathcal{H}_m$. We write $\mathcal{K}\varphi$ for the function

$$\mathcal{K}\varphi(\alpha) = \sum_{\beta \in \Omega_m} \mathcal{K}(\alpha, \beta) \varphi(\beta).$$

By the definition of $\mathcal{K}(\alpha, \beta)$, $\mathcal{K}\varphi$ does not depend on the order of the $m$ letters, for any $\varphi \in \mathcal{H}_m$. If $\bar{\mathcal{H}}_m$ denotes the linear subspace of functions $\varphi \in \mathcal{H}_m$ that are independent of the order of the $m$ letters, then any eigenvector of $\mathcal{K}$ must be in $\bar{\mathcal{H}}_m$. Clearly, constant functions give the eigenvalue $\lambda_0 = 1$.

We say that letter $j$ belongs to $\alpha \in \Omega_m$ (and write $j \in \alpha$) if the letter labeled $j$ appears in the configuration $\alpha$. For any $j = 1, \ldots, n$, we define $\chi_j$ to be the indicator function of the event $\{j \in \alpha\}$. To determine the spectrum, we may proceed as follows. Consider distinct letters $j_1, \ldots, j_k$, where $1 \le k \le m$. If at least one of the $k$ letters belongs to $\alpha$, then clearly $\mathcal{K}[\chi_{j_1} \cdots \chi_{j_k}](\alpha) = 0$. If none of them belong to $\alpha$, the probability of having the given $k$ letters in $\eta_{1,m}$, given that $\eta_{2,m} = \alpha$, is $\mathcal{K}[\chi_{j_1} \cdots \chi_{j_k}](\alpha) = (n-m-k)_m\,_{-k}/(n-m)_m$. Therefore, defining $\lambda_k$ as in (12), we have

(14) $$\mathcal{K}\chi_{j_1} \cdots \chi_{j_k} = |\lambda_k|(1 - \chi_{j_1}) \cdots (1 - \chi_{j_k}).$$

For $k \in \{1, \ldots, m\}$, let $\mathcal{U}_k$ denote the collection of functions of the form $\varphi = \chi_{j_1} \cdots \chi_{j_k}$, with $k$ distinct letters. Also, let $\mathcal{A}_k$ denote the linear span of

$$\{\mathcal{U}_0, \mathcal{U}_1, \ldots, \mathcal{U}_k\},$$

with $\mathcal{U}_0$ denoting the constant function $\varphi = 1$. Observe that, by (14), the subspaces $\mathcal{A}_k$ are invariant for $\mathcal{K}$, that is, $\mathcal{K}\mathcal{A}_k \subset \mathcal{A}_k$. Moreover, $\mathcal{A}_m$ coincides with $\bar{\mathcal{H}}_m$. Since $\mathcal{K}$ is self-adjoint w.r.t. $\langle \cdot, \cdot \rangle$, to establish (12) it is sufficient to prove that if $\varphi \in \mathcal{A}_k$ and $\varphi$ is orthogonal to $\mathcal{A}_{k-1}$, then

$$\langle \varphi, \mathcal{K}\varphi \rangle = \lambda_k \langle \varphi, \varphi \rangle. \tag{15}$$

To prove (15), observe that $\varphi \in \mathcal{A}_k$ must be of the form $\varphi = \phi + \psi$ with $\phi \in \mathcal{A}_{k-1}$ and $\psi$ a linear combination of functions in $\mathcal{U}_k$, say $\psi = \sum_i a_i \psi_i$. For each $\psi_i \in \mathcal{U}_k$, we use (14) to obtain $\mathcal{K}\psi_i = \widetilde{\phi}_i + \lambda_k \psi_i$, with $\widetilde{\phi}_i \in \mathcal{A}_{k-1}$. Setting $\widetilde{\phi} = \sum_i a_i \widetilde{\phi}_i$, we see that $\mathcal{K}\psi = \widetilde{\phi} + \lambda_k \psi$. If $\varphi$ is orthogonal to $\mathcal{A}_{k-1}$, we then have $\langle \varphi, \mathcal{K}\phi \rangle = \langle \varphi, \widetilde{\phi} \rangle = 0$. In particular,

$$\langle \varphi, \mathcal{K}\varphi \rangle = \lambda_k \langle \varphi, \psi \rangle = \lambda_k \langle \varphi, \varphi \rangle. \qquad \square$$

We now turn to a spectral gap estimate for an auxiliary block-dynamics. Suppose $m \leq n/N$ for some integer $N \geq 2$ and let $I_1, \ldots, I_N$ denote $N$ nonoverlapping $m$-blocks in $V_n$. Consider the Markov chain which picks uniformly at random one of the $N$ blocks and, given the configuration on that block, chooses uniformly at random a compatible configuration on the complement of that block. In symbols, we are looking at the symmetric transition matrix

$$\mathcal{P}(\sigma, \xi) = \frac{1}{N} \sum_{k=1}^{N} \nu[\eta = \xi | \eta = \sigma \text{ on } I_k], \qquad \sigma, \xi \in \Omega. \tag{16}$$

The next result is a natural extension of the estimates in Theorem 2.1 of [3] and Lemma 2.2 of [2].

PROPOSITION 2.2. *Let $\mu$ denote the lowest nonzero eigenvalue of $1 - \mathcal{P}$. Then,*

$$\mu \geq \frac{N-2}{N-1}. \tag{17}$$

PROOF. For all functions $f, g : \Omega \to \mathbb{R}$, we denote by $(f, g)$ the scalar product

$$(f, g) = \nu[fg] = \frac{1}{n!} \sum_{\eta \in \Omega} f(\eta) g(\eta). \tag{18}$$

Let $\pi_k : \Omega \to \Omega_m$ denote the projection $\eta \to \eta_{I_k}$, the restriction to the configuration on the block $I_k$. We define the subspace $\Gamma$ of sums of mean-zero functions of a single block:

$$\Gamma = \left\{ f = \sum_{k=1}^{N} g_k \circ \pi_k; g_k : \Omega_m \to \mathbb{R} \text{ with } \langle g_k, 1 \rangle = 0 \ \forall k \right\}, \tag{19}$$

with $\langle \cdot, \cdot \rangle$ as defined in (13). We write $\mathcal{P} = \frac{1}{N} \sum_{k=1}^{N} \mathcal{P}_k$, with
$$\mathcal{P}_k(\sigma, \xi) = \nu[\eta = \xi | \pi_k \eta = \pi_k \sigma].$$

In this way,
$$\mathcal{P}_k f(\sigma) = \sum_{\xi \in \Omega} \mathcal{P}_k(\sigma, \xi) f(\xi)$$

coincides with the conditional expectation of $f$, given the configuration in the block $I_k$, and is a function of the letters on the single block $I_k$ only. Therefore, $\mathcal{P} f \in \Gamma$ for every $f$ with $\nu[f] = 0$. Since $\mathcal{P}$ is self-adjoint w.r.t. $(\cdot, \cdot)$, $\mathcal{P} f = 0$ whenever $\nu[f] = 0$, and $f$ is orthogonal to $\Gamma$. To prove (17), it is then sufficient to establish

$$(20) \qquad (f, (1-\mathcal{P})f) \geq \frac{N-2}{N-1}(f,f), \qquad f \in \Gamma.$$

Given $f \in \Gamma$, $f = \sum_k g_k \circ \pi_k$, we define $\varphi_f = \sum_k g_k$, a function on $\Omega_m$. If $\mathcal{K}$ is defined by (11), we have
$$(g_j \circ \pi_j, g_k \circ \pi_k) = \begin{cases} \langle g_j, \mathcal{K} g_k \rangle, & k \neq j, \\ \langle g_k, g_k \rangle, & k = j. \end{cases}$$

We now compute
$$(f,f) = \sum_{k,j}(g_k \circ \pi_k, g_j \circ \pi_j)$$
$$= \sum_k \sum_{j: j \neq k} \langle g_k, \mathcal{K} g_j \rangle + \sum_k \langle g_k, g_k \rangle$$
$$(21)$$
$$= \sum_k \langle g_k, \mathcal{K} \varphi_f \rangle - \sum_k \langle g_k, \mathcal{K} g_k \rangle + \sum_k \langle g_k, g_k \rangle$$
$$= \langle \varphi_f, \mathcal{K} \varphi_f \rangle + \sum_k \langle g_k, (1-\mathcal{K}) g_k \rangle.$$

Similarly, observing that
$$\mathcal{P}_k(g_j \circ \pi_j) = \begin{cases} (\mathcal{K} g_j) \circ \pi_k, & k \neq j, \\ g_k \circ \pi_k, & k = j, \end{cases}$$

we compute, for every $k$,
$$(f, \mathcal{P}_k f) = \sum_{i,j}(g_i \circ \pi_i, \mathcal{P}_k(g_j \circ \pi_j))$$
$$= \sum_i \sum_{j: j \neq k}(g_i \circ \pi_i, (\mathcal{K} g_j) \circ \pi_k) + \sum_i (g_i \circ \pi_i, g_k \circ \pi_k)$$
$$(22)$$
$$= \sum_{i: i \neq k} \sum_{j: j \neq k} \langle g_i, \mathcal{K}^2 g_j \rangle + \sum_{j: j \neq k} \langle g_k, \mathcal{K} g_j \rangle + \sum_{i: i \neq k} \langle g_i, \mathcal{K} g_k \rangle + \langle g_k, g_k \rangle$$
$$= \langle \varphi_f, \mathcal{K}^2 \varphi_f \rangle + 2 \langle \varphi_f, \mathcal{K}(1-\mathcal{K}) g_k \rangle + \langle g_k, (1-\mathcal{K})^2 g_k \rangle.$$

Averaging over $k$ and rearranging terms, we finally have

$$(f, (1-\mathcal{P})f) = \frac{N-2}{N}\langle \varphi_f, \mathcal{K}(1-\mathcal{K})\varphi_f\rangle$$
(23)
$$+ \frac{1}{N}\sum_k \langle g_k, (1-\mathcal{K})((N-1)+\mathcal{K})g_k\rangle.$$

Note that (21) and (23) allow us to reduce the assertion (20) to spectral estimates involving only the operator $\mathcal{K}$. These, in turn, will follow from Lemma 2.1.

Consider the subspace $\mathcal{S} \subset \Gamma$ of symmetric functions,

(24) $$\mathcal{S} = \left\{ f = \sum_{k=1}^N g \circ \pi_k; g : \Omega_m \to \mathbb{R} \text{ with } \langle g, 1 \rangle = 0 \right\}.$$

Since $\mathcal{S}$ is invariant for $\mathcal{P}$, that is, $\mathcal{P}\mathcal{S} \subset \mathcal{S}$, to prove (20) we may consider separately the cases $f \in \mathcal{S}$ and $f \perp \mathcal{S}$. When $f \in \mathcal{S}$, that is, $f = \sum_k g \circ \pi_k$, we have $\varphi_f = Ng$, and, rearranging terms in (21) and (23), we obtain

(25) $$(f, f) = N(N-1)\left\langle g, \left[\mathcal{K} + \frac{1}{N-1}\right]g\right\rangle,$$

(26) $$(f, (1-\mathcal{P})f) = (N-1)^2 \left\langle g, [1-\mathcal{K}]\left[\mathcal{K} + \frac{1}{N-1}\right]g\right\rangle.$$

Since $m \leq n/2$, it is easily seen from (12) that $|\lambda_{k+1}| \leq |\lambda_k|$, for any $k$. Moreover, since $Nm \leq n$, we have $|\lambda_1| \leq 1/(N-1)$, so that $\mathcal{K} + \frac{1}{N-1}$ is nonnegative, and we may define its square root $[\mathcal{K}+\frac{1}{N-1}]^{1/2}$. Set $\tilde{g} = [\mathcal{K}+\frac{1}{N-1}]^{1/2}g$, so that $(f,f) = N(N-1)\langle \tilde{g}, \tilde{g}\rangle$. Observe that $\langle \tilde{g}, 1\rangle = \langle g, [\mathcal{K}+\frac{1}{N-1}]^{1/2}1\rangle = 0$ (since $[\mathcal{K}+\frac{1}{N-1}]^{1/2}$ is self-adjoint and $g$ is orthogonal to constants). Therefore,

$$\langle \tilde{g}, \mathcal{K}\tilde{g}\rangle \leq \lambda_2 \langle \tilde{g}, \tilde{g}\rangle,$$

since $\lambda_1$ is negative and $\lambda_2$ is the largest positive eigenvalue (other than $\lambda_0 = 1$).

From (25) and (26), we then have

(27) $$\begin{aligned}(f, (1-\mathcal{P})f) &= (N-1)^2\langle \tilde{g}, (1-\mathcal{K})\tilde{g}\rangle \\ &\geq (N-1)^2(1-\lambda_2)\langle \tilde{g}, \tilde{g}\rangle \\ &= \frac{N-1}{N}(1-\lambda_2)(f,f).\end{aligned}$$

From (12), we have $\lambda_2 = m(m-1)/(n-m)(n-m-1)$. In particular, for any $N \geq 2$ and $m \geq 1$ with $Nm \leq n$,

(28) $$\lambda_2 \leq \frac{1}{(N-1)^2}.$$

It follows that $\frac{N-1}{N}(1-\lambda_2) \geq \frac{N-2}{N-1}$. From (27), we conclude that

$$(f, (1-\mathcal{P})f) \geq \frac{N-2}{N-1}(f, f),$$

which proves the claim for $f \in \mathcal{S}$.

It remains to study the case $f \in \mathcal{S}^\perp$. Let $u$ be a generic element of $\mathcal{S}$, that is, $u = \sum_k u_0 \circ \pi_k$, for some $u_0$ with $\langle u_0, 1 \rangle = 0$. Computing as in (21), one has

$$(f, u) = (N-1)\left\langle \left[\mathcal{K} + \frac{1}{N-1}\right]\varphi_f, u_0 \right\rangle.$$

Since $u_0$ is an arbitrary mean-zero function, requiring $f \in \mathcal{S}^\perp$ [i.e. $(f, u) = 0$ for all $u \in \mathcal{S}$] implies that $[\mathcal{K} + \frac{1}{N-1}]\varphi_f$ is a constant. It follows that

$$\begin{aligned}
\left[\mathcal{K} + \frac{1}{N-1}\right]\varphi_f &= \left\langle \left[\mathcal{K} + \frac{1}{N-1}\right]\varphi_f, 1 \right\rangle \\
&= \left\langle \varphi_f, \left[\mathcal{K} + \frac{1}{N-1}\right]1 \right\rangle \\
&= \frac{N}{N-1}\langle \varphi_f, 1 \rangle = 0,
\end{aligned}$$

where we use the fact that $\langle \varphi_f, 1 \rangle = 0$ (recall that each $g_j$ has mean zero by assumption). In particular, this shows that

(29) $$\frac{N-2}{N}\langle \varphi_f, \mathcal{K}(1-\mathcal{K})\varphi_f \rangle = \frac{N-2}{N-1}\langle \varphi_f, \mathcal{K}\varphi_f \rangle.$$

The above identity says that the first term in the r.h.s. of (23) equals $(N-2)/(N-1)$ times the first term in the r.h.s. of (21).

Let us now look at the second term in the r.h.s. of (23). If we define the functions $\hat{g}_k = (1-\mathcal{K})^{1/2}g_k$, we have

(30) $$\frac{1}{N}\sum_k \langle g_k, (1-\mathcal{K})((N-1) + \mathcal{K})g_k \rangle = \frac{1}{N}\sum_k \langle \hat{g}_k, [(N-1) + \mathcal{K}]\hat{g}_k \rangle.$$

Using $\mathcal{K} \geq \lambda_1 \geq -1/(N-1)$, we obtain that (30) is estimated from below by

(31) $$\frac{N-2}{N-1}\sum_k \langle \hat{g}_k, \hat{g}_k \rangle = \frac{N-2}{N-1}\sum_k \langle g_k, (1-\mathcal{K})g_k \rangle.$$

This is precisely $(N-2)/(N-1)$ times the second term in the r.h.s. of (21). In conclusion, (29) and (31) show that

(32) $$(f, (1-\mathcal{P})f) \geq \frac{N-2}{N-1}(f, f), \qquad f \in \mathcal{S}^\perp.$$

From (27) and (32), we obtain the claim (17), and the proof is complete. $\square$

2.3. *The block-average dynamics.* Let $\ell$ and $N$ be given integers and suppose that we have exactly $n = N\ell$ vertices. We partition the vertex set $V_n$ by means of $N$ nonoverlapping $\ell$-blocks $I_{1,\ell}, \ldots, I_{N,\ell}$. There are $\ell$ choices for such partitioning, but, by symmetry, our estimates will not depend on the choice. For the sake of definiteness, we take the $\ell$-blocks $I_{k,\ell} = \{(k-1)\ell + 1, (k-1)\ell + 2, \ldots, k\ell\}$.

The block-average dynamics is the continuous time Markov chain obtained as follows. There is an independent rate-1 Poisson clock at each block. When block $I_{i,\ell}$ rings, a further block $I_{j,\ell}$ is chosen uniformly at random (with replacement, i.e. the choice may be $I_{i,\ell}$ itself). If $i = j$, we do nothing. If $i \neq j$, we choose uniformly the new configuration on $I_{i,\ell} \cup I_{j,\ell}$ among all configurations compatible with the letters outside the blocks $I_{i,\ell}$ and $I_{j,\ell}$.

The above defined process is reversible with respect to the uniform probability $\nu$, and its Dirichlet form can be written as

$$\mathcal{D}(f,f) = \frac{1}{N} \sum_{i=1}^{N} \sum_{j=1}^{N} \nu[(A_{i,j}f)^2], \qquad (33)$$

where the average gradient is given by

$$A_{i,j}f(\sigma) = \sum_{\xi \in \Omega} (f(\xi) - f(\sigma))\nu[\eta = \xi | \pi_k \eta = \pi_k \sigma \forall k \neq i, j], \qquad (34)$$

for $i \neq j$, and where we agree that $A_{i,i}f = 0$. Also, note that

$$\nu[(A_{i,j}f)^2] = \nu[\mathrm{Var}_{i,j}(f)],$$

where $\mathrm{Var}_{i,j}$ stands for the variance w.r.t. $\nu[\cdot | \pi_k \eta \forall k \neq i, j]$, the conditional probability obtained by freezing the configuration in all blocks $I_{k,\ell}$, $k \neq i, j$.

PROPOSITION 2.3. *Let $\gamma(N, \ell)$ denote the relaxation time of the block-average dynamics. Then,*

$$\gamma(N, \ell) = 1, \qquad (35)$$

*for every $\ell \geq 1$, $N \geq 2$, $n = N\ell$.*

PROOF. We first observe that $\gamma(2, \ell) = 1$ for all $\ell \geq 1$. Indeed, when $N = 2$, we have $A_{1,2}f = A_{2,1}f = \nu[f] - f$, so that $\mathcal{D}(f,f) = E[(f - \nu[f])^2] = \mathrm{Var}(f)$.

We turn to the case $N \geq 3$. Using the notation (18) of Proposition 2.2, we write, for any $f : \Omega \to \mathbb{R}$,

$$(f, (1 - \mathcal{P})f) = \frac{1}{N} \sum_{j=1}^{N} (f, f - \mathcal{P}_j f). \qquad (36)$$

Note that $\mathcal{P}_j f(\eta)$ coincides with the function $\eta \to \nu[f|\pi_j \eta]$, the conditional expectation of $f$, given the configuration on $I_{j,\ell}$. In particular, denoting by $\mathrm{Var}_j(f)$ the variance of $f$ w.r.t. $\nu[\cdot|\pi_j \eta]$, we have

$$(f, f - \mathcal{P}_j f) = (f - \mathcal{P}_j f, f - \mathcal{P}_j f) = \nu[\mathrm{Var}_j(f)].$$

Once the $\ell$ variables in $I_{j,\ell}$ are frozen, the measure $\nu[\cdot|\pi_j \eta]$ is the uniform measure on all $(n-\ell)!$ permutations of $n-\ell$ letters over $n-\ell$ vertices. For any $\eta \in \Omega$, we may therefore estimate, by definition of $\gamma(N,\ell)$,

$$\mathrm{Var}_j(f) \leq \frac{\gamma(N-1,\ell)}{N-1} \sum_{i \neq j} \sum_{k \neq j} \nu[(A_{i,k}f)^2|\pi_j \eta].$$

Taking $\nu$-expectation, we then have

$$(f, f - \mathcal{P}_j f) \leq \frac{\gamma(N-1,\ell)}{N-1} \sum_{i \neq j} \sum_{k \neq j} \nu[(A_{i,k}f)^2].$$

Averaging over $j$ and observing that

$$\sum_j \sum_{i \neq j} \sum_{k \neq j} \nu[(A_{i,k}f)^2] = N(N-2)\mathcal{D}(f,f),$$

we obtain

(37) $$\frac{1}{N}\sum_{j=1}^N (f, f - \mathcal{P}_j f) \leq \gamma(N-1,\ell)\frac{N-2}{N-1}\mathcal{D}(f,f).$$

If $(f,1) = 0$, we know by Proposition 2.2 that $(f,(1-\mathcal{P})f) \geq \frac{N-2}{N-1}(f,f)$. From (36) and (37), this implies

$$(f,f) \leq \gamma(N-1,\ell)\mathcal{D}(f,f).$$

In conclusion, we have shown that

(38) $$\gamma(N,\ell) \leq \gamma(N-1,\ell), \qquad N \geq 3.$$

This shows that $\gamma(N,\ell) \leq \gamma(2,\ell) = 1$ for any $N \geq 2$, $\ell \geq 1$.

To prove a lower bound on $\gamma(N,\ell)$, let $f$ denote the indicator function of the event that letter 1 belongs to block $I_{1,\ell}$. Clearly, the expectation of $f$ is $\frac{1}{N}$, and $\mathrm{Var}(f) = \frac{N-1}{N^2}$. On the other hand, one can easily check that, for every $j \neq 1$, one has $\nu[(A_{1,j}f)^2] = \frac{1}{2N}$. Therefore, for this function,

$$\mathcal{D}(f,f) = \frac{2}{N}\sum_{j:j \neq 1} \nu[(A_{1,j}f)^2] = \frac{N-1}{N^2} = \mathrm{Var}(f),$$

which implies $\gamma(N,\ell) \geq 1$. $\square$

2.4. *Extensions of the block-average dynamics.* Here we extend the dynamics defined by (33) to the case where $n$ is not a multiple of $\ell$. We thus consider the case $n = N\ell + m$, $1 \leq m \leq \ell - 1$. We say that a collection $D = \{I_j\}$ of subsets of $V_n$ is an $\ell$-partition if:

(i) $V_n = \bigcup_j I_j$ and $I_j \cap I_k = \varnothing$ for $j \neq k$,
(ii) all the $I_j$'s but one are $\ell$-blocks,
(iii) the remaining $I_j$ is an $m$-block.

Any $\ell$-partition is therefore made of $N$ $\ell$-blocks (for some integer $N$) and one $m$-block. An $\ell$-partition is said to be of type 1 if vertex $x = 1$ is the left endpoint of one of the blocks $I_j$. We use symbols $I_{j,\ell}$, $j = 1, \ldots, N$, for the $\ell$-blocks and $I_{N+1,m}$ for the $m$-block. Moreover, there are exactly $N+1$ $\ell$-partitions of type 1, depending on the position of $I_{N+1,m}$ relative to the $I_{j,\ell}$, and we denote by $D_k$ the $\ell$-partition of type 1 for which $I_{N+1,m}$ is the $k$th block in the partition starting from vertex 1. The elements of $D_k$ are denoted by $I_{j,\ell}^{(k)}$ and $I_{N+1,m}^{(k)}$ ($k = 1, \ldots, N+1$).

Let $\mathcal{D}_k(f,f)$ denote the Dirichlet form

$$\mathcal{D}_k(f,f) = \frac{1}{N} \sum_{i=1}^{N} \sum_{j=1}^{N} \nu[(A_{i,j}^{(k)} f)^2], \tag{39}$$

where the average gradients $A_{i,j}^{(k)}$ are defined by (34), with $I_{j,\ell}^{(k)}$ in place of $I_{j,\ell}$. Since the $m$-block $I_{N+1,m}^{(k)}$ is never updated, the Dirichlet form (39) is nonergodic. However, the average over $k$ of $\mathcal{D}_k$ is ergodic, and we have

PROPOSITION 2.4. *For any $\ell \geq 1$, $N \geq 2$ and $m \leq \ell - 1$ such that $n = N\ell + m$, for every function $f : \Omega \to \mathbb{R}$,*

$$\operatorname{Var}(f) \leq \tfrac{3}{2} \widehat{\mathcal{D}}(f,f), \tag{40}$$

*where*

$$\widehat{\mathcal{D}}(f,f) := \frac{1}{N+1} \sum_{k=1}^{N+1} \mathcal{D}_k(f,f). \tag{41}$$

PROOF. For $k = 1, \ldots, N+1$, let $J_k$ denote the $m$-block $I_{N+1,m}^{(k)}$ of the partition $D_k$. Also, let $\pi_k^*$ denote the projection $\eta \to \eta_{J_k}$. Note that, by construction, the $m$-blocks $J_k$ have no overlap. We may then apply Proposition 2.2, defining $\mathcal{P}$ as in (16), with $I_k := J_k$ and $N$ replaced by $N+1$ (for every $N \geq 2$, we have here $N+1$ blocks of length $m$). We then see that, for any $f$,

$$\operatorname{Var}(f) \leq \frac{N}{N-1}(f,(1-\mathcal{P})f) \leq \frac{3}{2}(f,(1-\mathcal{P})f). \tag{42}$$

We now observe, as in the proof of Proposition 2.3, that for each $k$, the expression $(f, (1-\mathcal{P}_k)f)$ can be written as $\nu[\mathrm{Var}_k(f)]$, where $\mathrm{Var}_k(f)$ stands for the variance of $f$ w.r.t. $\nu(\cdot|\pi_k^*\eta)$, the probability obtained by freezing the $m$-block $J_k$. Clearly, for each $\eta$, the measure $\nu(\cdot|\pi_k^*\eta)$ is uniform over the configurations of $n - m = N\ell$ given letters. For each $k$ and $\eta$, we may then apply Proposition 2.3 to the system of $N$ blocks $I_{j,\ell}^{(k)}$, $j = 1, \ldots, N$, and estimate

$$(43) \qquad \mathrm{Var}_k(f) \leq \frac{1}{N} \sum_{i=1}^{N} \sum_{j=1}^{N} \nu[(A_{i,j}^{(k)} f)^2 | \pi_k^*\eta].$$

Taking $\nu$-expectation and averaging over $k$, we therefore obtain

$$(44) \qquad (f, (1-\mathcal{P})f) \leq \frac{1}{N+1} \sum_{k=1}^{N+1} \mathcal{D}_k(f, f).$$

The proof is complete. $\square$

2.5. *Comparison estimates.* The first comparison allows us to progress from the block-average dynamics to an intermediate block dynamics, the moves of which are (1) averaging of adjacent blocks and (2) exchanging of adjacent blocks. The exchange moves are speeded up by a factor $N^2$. We refer to this intermediate process as the local average-exchange dynamics.

The setting is as in Proposition 2.4 above, with $\mathcal{D}_k$ defined by (39). Given the type 1 $\ell$-partition $D_k$, we agree to rename the $\ell$-blocks $I_{j,\ell}^{(k)}$ so that $I_{1,\ell}^{(k)}$ comes after the $m$-block $I_{N+1,m}^{(k)}$, $I_{2,\ell}^{(k)}$ after $I_{1,\ell}^{(k)}$ and so on. The exchange gradients $E_{i,i+1}^{(k)}$ appearing in the statement below are defined by

$$(45) \qquad E_{i,i+1}^{(k)} f(\eta) = f(\eta^{(i,i+1)}) - f(\eta),$$

where, for any $i, j$, $\eta^{(i,j)}$ denotes the configuration $\eta$ after the interchange of block $I_{i,\ell}^{(k)}$ with $I_{j,\ell}^{(k)}$, that is, $\pi_k \eta^{(i,j)} = \pi_k \eta$ for all $k \neq i, j$, while $\pi_j \eta^{(i,j)} = \pi_i \eta$ and $\pi_i \eta^{(i,j)} = \pi_j \eta$. We agree that $\eta^{(i,i)} = \eta$.

LEMMA 2.5. *For every function $f: \Omega \to \mathbb{R}$, for every $k = 1, \ldots, N+1$,*

$$(46) \qquad \mathcal{D}_k(f, f) \leq 3N^2 \sum_{i=1}^{N-1} \nu[(E_{i,i+1}^{(k)} f)^2] + \tfrac{1}{2} \sum_{i=1}^{N-1} \nu[(A_{i,i+1}^{(k)} f)^2].$$

PROOF. Let $T_{i,j}: \Omega \to \Omega$ denote the transformation $\eta \to \eta^{(i,j)}$. Then, assuming $i < j$, we see that

$$T_{i,j} = T_{i,i+1} T_{i+1,i+2} \cdots T_{j-2,j-1} T_{j-1,j} T_{j-2,j-1} \cdots T_{i+1,i+2} T_{i,i+1}.$$

Define $H = j - i$ and set
$$U_h = T_{i+h-1,i+h} \cdots T_{i+1,i+2} T_{i,i+1}, \qquad h = 1, \ldots, H.$$

Also, set
$$S_h = T_{j-h-1,j-h} \cdots T_{j-2,j-1} U_H, \qquad h = 1, \ldots, H-1.$$

In this way, setting $S_0 = U_H$ and $U_0 = 1$, we write
$$\begin{aligned}
E_{i,j}^{(k)} f(\eta) &= f(\eta^{(i,j)}) - f(\eta) \\
&= f(S_{H-1}\eta) - f(U_H\eta) + f(U_H\eta) - f(\eta) \\
&= \sum_{h=1}^{H-1} [f(T_{j-h-1,j-h} S_{h-1}\eta) - f(S_{h-1}\eta)] \\
&\quad + \sum_{h=1}^{H} [f(T_{i+h-1,i+h} U_{h-1}\eta) - f(U_{h-1}\eta)] \\
&= A(f) + B(f), \qquad \text{say}.
\end{aligned}$$

We then estimate
$$\nu[(E_{i,j}^{(k)} f)^2] \le 2\nu[(A(f))^2] + 2\nu[(B(f))^2].$$

By Schwarz's inequality,
$$\begin{aligned}
\nu[(A(f))^2] &\le (H-1) \sum_{h=1}^{H-1} \nu[(f \circ T_{j-h-1,j-h} S_{h-1} - f \circ S_{h-1})^2] \\
&= (H-1) \sum_{h=1}^{H-1} \nu[(f \circ T_{j-h-1,j-h} - f)^2] \\
&\le N \sum_{i=1}^{N-1} \nu[(E_{i,i+1}^{(k)} f)^2],
\end{aligned}$$

where we have used the invariance of $\nu$ under all transformations $S_h$. The same reasoning gives the same estimate for $\nu[(B(f))^2]$. This shows that, for each couple $i,j$, we obtain the following estimate for the exchange terms:

(47) $$\nu[(E_{i,j}^{(k)} f)^2] \le 4N \sum_{i=1}^{N-1} \nu[(E_{i,i+1}^{(k)} f)^2].$$

We now have to estimate the average terms. To this end, observe that, if the symbol $\omega_{i,j}(\eta)$ is used to denote the configuration $\eta$ outside the two blocks $I_{i,\ell}^{(k)}$ and $I_{j,\ell}^{(k)}$, we may regard $A_{i,j}^{(k)} f$ as the function $A_{i,j}^{(k)} f(\eta) = \nu[f|\omega_{i,j}(\eta)] -$

$f$. Now a simple change of variables shows that for any other block, $I_{h,\ell}^{(k)}$, $h \neq i, j$,

$$\nu[f|\omega_{i,j}(\eta)] = \nu[f^{(j,h)}|\omega_{i,h}(\eta^{(j,h)})],$$

where, as usual, $\eta^{(j,h)}$ stands for the configuration in which the blocks $I_{j,\ell}^{(k)}, I_{h,\ell}^{(k)}$ have been exchanged and $f^{(j,h)}(\eta) = f(\eta^{(j,h)})$. We then have, using the invariance of $\nu$ under the exchange $\eta \to \eta^{(j,h)}$,

$$\nu[(A_{i,j}^{(k)} f)^2] = \nu[(\nu[f^{(j,h)}|\omega_{i,h}(\eta^{(j,h)})] - f)^2]$$
$$= \nu[(\nu[f^{(j,h)}|\omega_{i,h}(\eta)] - f^{(j,h)})^2]$$
$$= \nu[(A_{i,h}^{(k)} f^{(j,h)})^2].$$

Moreover,

$$\nu[(A_{i,h}^{(k)} f^{(j,h)})^2] \leq 2\nu[(A_{i,h}^{(k)} f)^2] + 2\nu[(A_{i,h}^{(k)} (f^{(j,h)} - f))^2],$$

and Jensen's inequality implies

$$\nu[(A_{i,h}^{(k)} (f^{(j,h)} - f))^2] \leq \nu[(f^{(j,h)} - f)^2].$$

The last expression is nothing but $\nu[(E_{j,h}^{(k)} f)^2]$, and we can estimate it by (47). Summarizing, choosing, for example, $h = j + 1$ we have obtained

$$(48) \qquad \nu[(A_{i,j}^{(k)} f)^2] \leq 2\nu[(A_{j,j+1}^{(k)} f)^2] + 8N \sum_{l=1}^{N-1} \nu[(E_{l,l+1}^{(k)} f)^2].$$

Recalling the expression (39) for the Dirichlet form $\mathcal{D}_k(f,f)$, the bounds (47) and (48) imply the desired assertion (46). $\square$

The next task is to rewrite the exchange and average contributions to the local block dynamics appearing in the r.h.s. of (46) in terms of the inversions $\eta \to \eta^{x,\ell}$. We start by writing each exchange term in terms of three inversions. To this end, we denote by $z$ the left end vertex of $I_{i,\ell}^{(k)}$, so that $z + \ell$ is the left end vertex of $I_{i+1,\ell}^{(k)}$.

LEMMA 2.6. *For every function $f : \Omega \to \mathbb{R}$, for every $k = 1, \ldots, N+1$,*

$$(49) \quad \nu[(E_{i,i+1}^{(k)} f)^2] \leq 3\nu[(R_{z,2\ell-1} f)^2 + (R_{z+\ell,\ell-1} f)^2 + (R_{z,\ell-1} f)^2].$$

PROOF. Denote by $W_{x,h} : \Omega \to \Omega$ the inversion transformation $\eta \to \eta^{x,\ell}$. It is clear that $W_{z,2\ell-1} W_{z+\ell,\ell-1} W_{z,\ell-1} \eta$ coincides with $\eta^{(i,i+1)}$. Therefore,

$$E_{i,i+1}^{(k)} f(\eta) = R_{z,2\ell-1} f(W_{z+\ell,\ell-1} W_{z,\ell-1} \eta)$$
$$+ R_{z+\ell,\ell-1} f(W_{z,\ell-1} \eta) + R_{z,\ell-1} f(\eta).$$

The conclusion now follows from Schwarz's inequality and the invariance of $\nu$ under all transformations $W_{x,h}$. □

Our last comparison shows how to bound each average term in (46) by means of inversions.

LEMMA 2.7. *For every function* $f:\Omega \to \mathbb{R}$, *for every* $k = 1, \ldots, N+1$,

$$(50) \quad \nu[(A_{i,i+1}^{(k)}f)^2] \leq \frac{1}{2\ell} \sum_{x \in I_i^{(k)} \cup I_{i+1}^{(k)}} \sum_{l \leq 2\ell} \nu[(R_{x,l}f)^2 + (R_{x+1,l-1}f)^2].$$

PROOF. The starting point is the following standard spectral gap estimate for random transpositions. Consider a system with $d$ vertices and $d$ letters. Let $\nabla_{x,y}$ denote the gradient associated to the transposition of the letters at vertices $x$ and $y$, that is, $\nabla_{x,y}f(\eta) = f(\eta^{x,y}) - f(\eta)$, with $\eta^{x,y}$ denoting the configuration identical to $\eta$ out of $x,y$ and such that $(\eta^{x,y})_x = \eta_y$ and $(\eta^{x,y})_y = \eta_x$. It is then well known (see, e.g., Theorem 5.1 in [2] for a simple proof) that, for all $d \geq 2$ and all functions $f$,

$$(51) \quad \mathrm{Var}_d(f) \leq \frac{1}{4d} \sum_{x,y} \nu_d[(\nabla_{x,y}f)^2],$$

where $\nu_d$ is the uniform measure over the $d!$ permutations, and $\mathrm{Var}_d$ denotes the variance w.r.t. $\nu_d$.

We want to apply this bound to the system with $d = 2\ell$ which is obtained by freezing all the letters in the complement of $I_i^{(k)} \cup I_{i+1}^{(k)}$. Recall the notation $A_{i,i+1}^{(k)}f(\eta) = \nu[f \mid \omega_{i,i+1}(\eta)] - f$ that was used in the proof of Lemma 2.5. We then have, for each $\eta$ and every $f: \Omega \to \mathbb{R}$,

$$(52) \quad \mathrm{Var}(f \mid \omega_{i,i+1}(\eta)) \leq \frac{1}{8\ell} \sum_{x,y \in I_i^{(k)} \cup I_{i+1}^{(k)}} \nu[(\nabla_{x,y}f)^2 \mid \omega_{i,i+1}(\eta)].$$

Clearly, each transposition is written as the composition of two inversions: if, for example, $y = x+h$, then $\nabla_{x,x+h}f(\eta) = R_{x,h}f(\eta^{x+1,x+h-1}) + R_{x+1,h-1}f(\eta)$. Recalling that $\nu[(A_{i,i+1}^{(k)}f(\eta))^2] = \nu[\mathrm{Var}(f \mid \omega_{i,i+1}(\eta))]$, the bound (50) is seen to be a straightforward consequence of (52). □

2.6. *Proof of the upper bound in Theorem* 1.1. We now combine the spectral gap estimate of Proposition 2.4 with the three comparison lemmas above. Recall that we are looking for a bound of the form $\tau(n,L) \leq C(n \vee \frac{n^3}{L^3})$, that is,

$$(53) \quad \mathrm{Var}(f) \leq C\left(n \vee \frac{n^3}{L^3}\right) \mathcal{E}(f,f),$$

for an arbitrary function $f$, with $\mathcal{E}(f,f)$ given by

$$\mathcal{E}(f,f) = \frac{1}{2}\frac{1}{nL}\sum_{x\in V_n}\sum_{\ell\leq L}\nu[(R_{x,\ell}f)^2]. \tag{54}$$

For every $1\leq \ell \leq L$, we use Proposition 2.4 and Lemma 2.5 to write

$$\mathrm{Var}(f) \leq CN^2\mathrm{Ex}(f) + C\mathrm{Av}(f), \tag{55}$$

with the exchange and average terms given by

$$\mathrm{Ex}(f) := \frac{1}{N+1}\sum_{k=1}^{N+1}\sum_{i=1}^{N-1}\nu[(E_{i,i+1}^{(k)}f)^2], \tag{56}$$

$$\mathrm{Av}(f) := \frac{1}{N+1}\sum_{k=1}^{N+1}\sum_{i=1}^{N-1}\nu[(A_{i,i+1}^{(k)}f)^2]. \tag{57}$$

We start by estimating the average term $\mathrm{Av}(f)$. Using Lemma 2.7 and observing that for every $x \in V_n$, on any given partition $D_k$ the number of indices $i$ such that $x \in I_i^{(k)} \cup I_{i+1}^{(k)}$ is at most 2, we have

$$\mathrm{Av}(f) \leq \frac{2}{\ell}\sum_{x\in V_n}\sum_{l\leq 2\ell}\nu[(R_{x,l}f)^2 + (R_{x+1,l-1}f)^2]$$

$$\leq \frac{4}{\ell}\sum_{x\in V_n}\sum_{l\leq 2\ell}\nu[(R_{x,l}f)^2].$$

Therefore, for any $\delta > 0$ and $\ell \leq n$ such that $\delta L \leq \ell \leq \frac{1}{2}L$, we have

$$\mathrm{Av}(f) \leq 8\delta^{-1}n\mathcal{E}(f,f). \tag{58}$$

Summarizing, taking, for example, $\delta = 0.1$ and adjusting the value of the constant $C$ in (55), we have obtained that if $\delta L \leq \ell \leq \frac{1}{2}L$, then, for every $f$,

$$\mathrm{Var}(f) \leq CN^2\mathrm{Ex}(f) + Cn\mathcal{E}(f,f). \tag{59}$$

We now turn to an estimate of the exchange terms (56). Observe that, so far, only partitions of type 1 enter the definition of $\mathrm{Ex}(f)$; see (56). However, by symmetry, we could just as well use, for any $y \in V_n$, partitions of type $y$ defined as those $\ell$-partitions $\{I_j\}$ of $V_n$ for which the vertex $y$ is the left endpoint of one of the blocks $I_j$. Again, there are exactly $N+1$ $\ell$-partitions of type $y$, say $D_{k,y}$, $k=1,\ldots,N+1$. Let us denote by $I_{j,\ell}^{(k,y)}$ the corresponding $\ell$-blocks and define the associated exchange gradients $E_{i,j}^{(k,y)}$ as in (45) and (34) with $I_{j,\ell}^{(k,y)}$ in place of $I_{j,\ell}$. We may then replace $\mathrm{Ex}(f)$ in (59) by the following average over partitions of type $y$, for $y=1,\ldots,\ell$:

$$\frac{1}{\ell}\sum_{y=1}^{\ell}\mathrm{Ex}_y(f), \qquad \mathrm{Ex}_y(f) := \frac{1}{N+1}\sum_{k=1}^{N+1}\sum_{i=1}^{N-1}\nu[(E_{i,i+1}^{(k,y)}f)^2]. \tag{60}$$

Now, each term $\text{Ex}_y(f)$ is estimated by means of Lemma 2.6. Moreover, it is easy to check that for every $x \in V_n$ and $k = 1, \ldots, N+1$, there is only one vertex $y \in \{1, \ldots, \ell\}$ such that $x$ is the left endpoint of one of the $I_{j,\ell}^{(k,y)}$ ($j = 1, \ldots, N$). It follows that

$$\frac{1}{\ell}\sum_{y=1}^{\ell} \text{Ex}_y(f) \leq \frac{3}{\ell} \sum_{x \in V_n} \nu[(R_{x,2\ell-1}f)^2 + (R_{x+\ell,\ell-1}f)^2 + (R_{x,\ell-1}f)^2].$$

We want to average this expression over $\delta L \leq \ell \leq \frac{1}{2}L$. To this end, choose $\delta = 0.1$ and $L_0$ such that the number of integers $\ell$ satisfying $\delta L \leq \ell \leq \frac{1}{2}L$ is greater than $L/4$ for every $L \geq L_0$. If $L \geq L_0$, we can therefore estimate

$$\frac{1}{\#\{\ell : \delta L \leq \ell \leq (1/2)L\}} \sum_{\delta L \leq \ell \leq (1/2)L} \frac{1}{\ell}\sum_{y=1}^{\ell} \text{Ex}_y(f) \leq C\frac{n}{L}\mathcal{E}(f,f),$$

for some universal constant $C$. Going back to (59), we first replace $\text{Ex}(f)$ by (60) and then take the average over $\delta L \leq \ell \leq \frac{1}{2}L$. In this way, recalling that $N = [n/\ell] \leq n/(\delta L)$ we obtain, for $L \geq L_0$ and a new universal constant $C$,

(61) $$\text{Var}(f) \leq C\frac{n^3}{L^3}\mathcal{E}(f,f) + Cn\mathcal{E}(f,f) \leq 2C\left(n \vee \frac{n^3}{L^3}\right)\mathcal{E}(f,f).$$

It remains to consider the case $L \leq L_0$. In this case, it is sufficient to prove that $\tau(n, L) \leq Cn^3$, for some universal constant $C$. Recall that, when $L = 1$, $\tau(n, 1) \leq Cn^3$, by the well-known diffusive bound for local transpositions (see, e.g., [2]). On the other hand, the obvious bound

$$\mathcal{E}(f,f) \geq \frac{1}{nL_0} \sum_{x \in V_n} \nu[(R_{x,1}f)^2]$$

shows that $\tau(n, L) \leq L_0 \tau(n, 1)$, for any $L \leq L_0$. This concludes the proof of the upper bound in Theorem 1.1.

2.7. *Proof of the upper bound in Theorem* 1.2. This is a straightforward consequence of the previous bounds on $\tau(n, L)$. Indeed, there is some $C_1 < \infty$ such that $\theta^{\ell-1} \geq 1/C_1$ for any $\ell \leq (1-\theta)^{-1}$, $\theta \in (0,1)$. Therefore, choosing an integer $L \leq (1-\theta)^{-1}$ with $L^{-1} \leq C_2(1-\theta)$ and removing all inversions with $\ell > L$, we have

$$\tau(n,\theta) \leq C_1 C_2 \tau(n,L) \leq C_1 C_2^4 (n \vee [n(1-\theta)]^3).$$

**3. Lower bounds.** Here we are going to find the appropriate slow modes of relaxation. As already discussed in [6, 8], the bound $O(n \vee \frac{n^3}{L^3})$ corresponds to the competition of two different phenomena: the diffusive scale $O(\frac{n^3}{L^3})$ is related to the transport of local information, while the scale $O(n)$ originates from the difficulty of separating two adjacent letters during the time evolution (described in terms of conserved edges in [6]). For the sake of completeness, below we explicitly derive these bounds, although some of them appear already (in some cases with explicit constants) in [6].

3.1. *Proof of the lower bound in Theorem* 1.1. Let the $n$ letters be labeled by the integers from 1 to $n$ and denote by $\xi_x$ the indicator function of the event {letter $n$ is at vertex $x$}. Let $g:[0,1] \to \mathbb{R}$ be a smooth function with $\int_0^1 g(t)\,dt = 0$, $\int_0^1 g(t)^2\,dt = 1$. Set

(62) $$\psi(\eta) = \sum_{x \in V_n} g(x/n)\xi_x.$$

As $n \to \infty$, we have $\nu[\psi] \to 0$ and $\nu[\psi^2] \to 1$. We may therefore estimate $\mathrm{Var}(\psi) \geq \frac{1}{2}$, if $n$ is large enough. The lower estimate $\tau(n,L) \geq \frac{n^3}{CL^3}$ is then a consequence of

LEMMA 3.1. *There exists $C < \infty$ such that*

(63) $$\mathcal{E}(\psi,\psi) \leq C \frac{L^3}{n^3} \int_0^1 g'(t)^2\,dt.$$

PROOF. Write $I_{x,\ell}$ for the $(\ell+1)$-block with left end $x$ and right end $x+\ell$. Observe that, for any $x \in V_n$ and $\ell \leq L$, we have

$$\sum_{y \in I_{x,\ell}} (\xi_y^{x,\ell} - \xi_y) = 0.$$

It follows that

$$R_{x,\ell}\psi = \sum_{y \in I_{x,\ell}} (g(y/n) - g(x/n))(\xi_y^{x,\ell} - \xi_y).$$

We now expand $g(y/n) - g(x/n) = g'(x/n)(y-x)/n + O(\ell^2/n^2)$. Observing that

$$\sum_{y \in I_{x,\ell}} \nu[(\xi_y^{x,\ell} - \xi_y)^2] \leq 2(\ell+1)/n,$$

and ignoring higher-order terms, we have

$$\nu[(R_{x,\ell}\psi)^2] \leq C g'(x/n)^2 \frac{\ell^3}{n^3}.$$

The proof is then completed by averaging over $x \in V_n$ and $\ell \leq L$. □

To prove the bound $\tau(n, L) \geq n/C$, we may proceed as follows. Let $\chi$ denote the indicator function of the event {letter 1 is adjacent to letter 2}. It is easily seen that, for every $n \geq 3$, the probability of this event is $\nu[\chi] = 2/(n-1)$. The variance is given by $\text{Var}(\chi) = 4(n-2)/(n-1)^2$. The desired estimate then follows from

LEMMA 3.2. *For every $n \geq 3$,*

$$\mathcal{E}(\chi, \chi) \leq \frac{16}{n(n-1)}. \tag{64}$$

PROOF. Observe that $(R_{x,\ell}\chi)^2 = \chi(1 - \chi^{x,\ell}) + \chi^{x,\ell}(1 - \chi)$. Let us analyze the contribution of $\chi(1 - \chi^{x,\ell})$ (the other term contributes the same by symmetry). If $1, 2$ are adjacent in $\eta$, but not in $\eta^{x,\ell}$, this implies that either $\eta_x \in \{1, 2\}$ or $\eta_{x+\ell} \in \{1, 2\}$. Therefore, for every $\ell$,

$$\sum_{x \in V_n} \chi(1 - \chi^{x,\ell}) \leq 4\chi.$$

Taking expectations and dividing by $n$ then gives the desired result. □

3.2. *Proof of the lower bound in Theorem* 1.1. For the bound $\tau(n, \theta) \geq [n(1-\theta)]^3/C$, we simply follow the computation in Lemma 3.1 to obtain

$$\mathcal{E}_\theta(\psi, \psi) \leq \frac{C(1-\theta)}{n^3} \sum_{\ell \leq n} \theta^{\ell-1} \ell^3 \int_0^1 g'(t)^2 \, dt.$$

The required bound now follows from the fact that $(1-\theta)^4 \sum_{\ell \leq n} \theta^{\ell-1} \ell^3$ is uniformly bounded in $\theta \in (0, 1)$. Finally, the bound $\tau(n, \theta) \geq n/C$ is straightforward, since Lemma 3.2 still holds when $\mathcal{E}(\chi, \chi)$ is replaced by $\mathcal{E}_\theta(\chi, \chi)$.

3.3. *A lower bound on the log-Sobolev constant.* Here we prove the estimate (10). Let $n$ be even and denote by $\xi$ the indicator function of the event

{the first $n/2$ letters occupy the first $n/2$ vertices}.

Then, $\nu[\xi] = \binom{n}{n/2}^{-1}$ and

$$\text{Ent}(\xi) = -\nu[\xi] \log \nu[\xi] \geq \frac{n}{C} \nu[\xi]. \tag{65}$$

On the other hand, only inversions $\eta \to \eta^{x,\ell}$ with $x \in \{n - L, \ldots, n\} \cup \{n/2 - L, \ldots, n/2\}$ can affect the value of $\xi$. In such cases, we have $\nu[(R_{x,\ell}\xi)^2] \leq 2\nu[\xi]$ for any $\ell \leq L$. It follows that

$$\mathcal{E}(\sqrt{\xi}, \sqrt{\xi}) = \mathcal{E}(\xi, \xi) \leq \frac{CL}{n} \nu[\xi]. \tag{66}$$

Combining (65) and (66), we arrive at the desired estimate.

N. CANCRINI
DIPARTIMENTO DI MATEMATICA PURE ED APPLICATA
UNIVERSITA DEGLI STUDI DI L'AQUILA
I-67100 L'AQUILA
ITALY
E-MAIL: nicoletta.cancrini@roma1.infn.it

P. CAPUTO
F. MARTINELLI
DIPARTIMENTO DI MATEMATICA
UNIVERSITA ROMA TRE
LARGO S. MURIALDO 1
00146 ROMA
ITALY
E-MAIL: caputo@mat.uniroma3.it
        martin@mat.uniroma3.it